\documentclass{epl}
\usepackage{amsmath} \usepackage{amsfonts} \usepackage{amssymb}
\usepackage{graphics,color}

\newtheorem{thm}{Theorem}

\title{Rigorous nonperturbative Ornstein-Zernike
theory for Ising ferromagnets}
\shorttitle{Nonperturbative OZ theory}
\author{M. Campanino\inst{1} \and D. Ioffe\inst{2} \and Y. Velenik\inst{3}}
\shortauthor{M. Campanino \etal}
\institute{
  \inst{1} Dipartimento di Matematica, Universit\`a di Bologna, piazza di
Porta S. Donato 5, I-40126  Bologna, Italy\thanks{E-mail: \email{campanin@dm.unibo.it}}\\
  \inst{2} Faculty of Industrial Engineering, Technion, Haifa 3200,
  Israel\thanks{E-mail: \email{ieioffe@ie.technion.ac.il}}\\
  \inst{3} UMR-CNRS 6632, CMI, 39 rue Joliot Curie, 13453 Marseille,
  France\thanks{E-mail: \email{velenik@cmi.univ-mrs.fr}}
}

\begin{document}

\maketitle

\begin{abstract}
We rigorously derive the Ornstein-Zernike asymptotics of the pair-correlation
functions for finite-range Ising ferromagnets in any dimensions and at any
temperature above critical.
\end{abstract}

The celebrated heuristic argument by Ornstein and Zernike~\cite{OrZe1914}
 implies  that the asymptotic form of the truncated two-point
density correlation function of simple fluids away from the critical region is given by
\begin{equation}
\label{eq:OZ}
G(\vec{r}) \simeq\frac{A_\beta}{\sqrt{|\vec{r}|^{d-1}}} \,e^{-
|\vec{r}|\xi_\beta },
\end{equation}
where the value of the inverse correlation length $\xi_\beta$
depends only on the density $\rho$, the inverse
temperature $\beta$ and the spatial dimension $d$. The original OZ
 approach  hinges on the assumption that the so called direct correlation
function $C(\cdot )$,  which is \textit{de facto} introduced
through the renewal type relation
\begin{equation}
\label{eq:OZrenewal}
G(\vec{r} ) = C(\vec{r} ) + \rho \int_{\mathbb{R}^d}
C( \vec{r}-\vec{r}_1 ) G(\vec{r}_1 )\mathrm{d}\vec{r}_1 ,
\end{equation}
is of an appropriately short range.

Because of the physical significance of both the conclusions and of
the underlying heuristic assumptions a number of works
(see e.g. ~\cite{AbKu1977,Pa1978,Si1981,BrFr1985,MiZh1996})  were devoted to
attempts to put the theory on a rigorous basis, that is to derive
\eqref{eq:OZ} directly from the microscopic
picture of intermolecular interactions.
Most of these works, however, were based on expansion/perturbation
techniques and required technical low density or
high/low temperature assumptions and, thereby, addressed the situation
when the parameters are far away from the critical region. Since 
the OZ theory was, above all, intended to describe the phenomenon 
of critical scattering it would be of interest to devise such rigorous approach
to \eqref{eq:OZ} which would rely  only on qualitative features of 
noncriticality such as, for example,  finite compressibility 
(or finite susceptibility
 in the context of ferromagnetic lattice models).

In this Letter, we present a fully nonperturbative derivation  
of the direction dependent 
analog of ~\eqref{eq:OZ} for finite-range ferromagnetic Ising model above
the critical temperature in any dimension. Let 
${\mathbf J}=\{J_v\}_{v\in\mathbb{Z}^d}$, be a collection of 
nonnegative real numbers  such that $J_v=J_{-v}$ and $J_v=0$ if $|v|>R$, where
$R$ is some finite number. The (formal) Hamiltonian is then of the form
\begin{equation}
\label{eq:Hamiltonian}
H(\sigma) = -\tfrac12\sum_{x,y \in \mathbb{Z}^d} J_{y-x} \sigma_x\sigma_y\,.
\end{equation}
Our approach pertains to the  high temperature region $\beta <\beta_c = \beta_c ({\mathbf J},d)$, which is the set of all $\beta$ such that the susceptibility 
$
\chi_\beta = 
\sum_{x\in\mathbb{Z}^d}  \langle \sigma_0 \sigma_x\rangle_\beta
$
 is finite. By the Simon-Lieb argument $\chi_\beta <\infty $ implies 
strict exponential decay of the two-point function. Equivalently, the
series
\begin{equation}
\label{eq:Suscept}
\chi_\beta (t) = 
\sum_{x\in\mathbb{Z}^d} \, e^{(t,x)} 
\langle \sigma_0 \sigma_x\rangle_\beta \, ,
\end{equation}
where $(\cdot ,\cdot )$ is the usual scalar product in ${\mathbb R}^d$, 
has a nonempty domain of convergence for each subcritical value 
of the inverse temperature $\beta <\beta_c$.
Note that, by an important result  
of Aizenman, Barsky and
Fern\'andez~\cite{AiBaFe1987}, $\beta_c$ is actually the usual critical temperature: The spontaneous magnetization is positive whenever $\beta >\beta_c$. 

In the sequel we shall use ${\mathbf K}_\beta$ to denote
the domain of convergence of \eqref{eq:Suscept}. From a purely geometric 
point of view 
 the direction dependent inverse correlation length $\xi_\beta$ 
is the support function of ${\mathbf K}_\beta$. In particular, 
 the dependence of $\xi_\beta (n)$ on the direction  
$n\in \mathbb{S}^{d-1}$ is encoded
in the geometry of $\partial {\mathbf K}_\beta$. 
\begin{thm}
\label{theorem} In any dimension $d\geq 1$ and 
for any ferromagnetic model \eqref{eq:Hamiltonian} the 
asymptotic decay of the two-point correlation function in the 
high-temperature region $\beta <\beta_c ({\mathbf J} ,d)$  is given 
by
\begin{equation}
\label{eq:OZising}
\langle\sigma_0\,\sigma_x\rangle_\beta \simeq
\frac{\Psi_\beta(n_x)}{\sqrt{|x|^{(d-1)}}}\;
e^{-\xi_\beta(n_x)\, |x|}\; (1+o(1))\,,
\end{equation}
where $n_x=x/|x|\in \mathbb{S}^{d-1}$ is the unit  vector in the direction 
of $x$, and the function $\Psi_\beta$ is strictly positive and
analytic.  Moreover,  
  the 
inverse correlation length $\xi_\beta (n)$ is an analytic function
of the direction  
$n\in \mathbb{S}^{d-1}$ in the sense that the boundary
$\partial \mathbf{K}_\beta$ 
of $\mathbf{K}_\beta$
is  locally 
analytic and strictly convex. Furthermore, the  
Gaussian curvature
$\kappa_\beta $
 of  $\partial\mathbf{K}_\beta$ is uniformly positive,
\begin{equation*}
\bar{\kappa}_\beta= \min_{t\in\partial\mathbf{K}_\beta}\kappa_\beta (t)~>~0.
\end{equation*}
\end{thm}

Alternatively, strict convexity and
analyticity  of the direction dependent inverse correlation length
$\xi_\beta (\cdot )$ could be formulated in terms  of
the geometry of  the unit sphere $\partial \mathbf{U}_\beta$ in the
$\xi_\beta$-norm.

A full proof of Theorem~\ref{theorem} can be found in~\cite{CaIoVe2002}.

In principle our approach pertains to any model in which the pair-correlation
function $g(x)$
admits a suitable graphical
 representation of the type
\begin{equation}
\label{eq:GR}
g(x) = \sum_{\lambda:\, 0\to x} q(\lambda)\,,
\end{equation}
where the sum runs over a family of paths or, more generally, pathlike objects
 connecting $0$ and $x$, possibly with
compatibility constraints, e.g., some
form of self-avoidance; let us call such paths
admissible. The weights $q(\,\cdot\,)$  are supposed to
be strictly positive and to possess  a variation of
the following four properties:
\begin{itemize}
\item \textbf{Strict exponential decay of the two-point function} There 
exists $C_1<\infty$ such that, for
all $x\in\mathbb{Z}^d \setminus \{0\}$,
\begin{equation}
\label{eq:PUB}
g(x)\,  =\,  
\sum_{\lambda:\, 0\to x} q(\lambda) \leq C_1\, e^{-\xi(x)}\,,
\end{equation}
where $\xi(x) = -\lim_{k\to\infty}(k|x|)^{-1}\, \log g([kx])$ is the inverse
correlation length.
\item\textbf{Finite energy condition:} For any pair of compatible 
paths $\lambda$ and $\eta$ define the conditional weight
$$
q(\lambda \,|\, \eta) = q(\lambda \amalg \eta ) / q(\eta)\,
$$
where $\lambda\amalg\eta$ denotes the concatenation of $\lambda$ and $\eta$. 
Then there exists a universal 
 finite constant $C_2<\infty$ such that the conditional
weights are controlled in terms of path sizes $|\lambda |$ as:
\begin{equation}
\label{eq:FE}
q(\lambda \,|\, \eta)\,  \geq \, e^{-C_2 |\lambda |}\, . 
\end{equation}
\item \textbf{Splitting property:} There exists $C_3<\infty$, such 
that, for all
$x,y\in\mathbb{Z}^d \setminus \{0\}$ with $x\neq y$,
\begin{equation}
\label{eq:Spl}
\sum_{\lambda:\, 0\to x\to y} q(\lambda) \leq C_2\, \sum_{\lambda:\, 0\to x}
q(\lambda)\; \sum_{\lambda:\, x\to y} q(\lambda)\,.
\end{equation}
\item \textbf{Exponential mixing :}
There exists $C_4<\infty$ and $\theta\in (0,1)$ such that, 
for any four paths $\lambda$, $\eta$,
$\gamma_1$ and $\gamma_2$, with $\lambda \amalg \eta \amalg \gamma_1$ and
$\lambda \amalg \eta \amalg \gamma_2$ both admissible,
\begin{equation}
\label{eq:ExD}
\frac{q(\lambda \,|\, \eta\amalg\gamma_1)}{q(\lambda \,|\, \eta\amalg\gamma_2)}
\leq \exp\{C_4\, \sum_{\substack{x\in \lambda\\y\in\gamma_1\cup\gamma_2}}
\theta^{|x-y|} \}\,.
\end{equation}
\end{itemize}
Many models enjoy a graphical representation
 of correlation functions of the
form~\eqref{eq:GR}. In perturbative regimes, cluster expansions
provide a generic  example. Nonperturbative 
examples include the random-cluster
representation for Potts (and other) models~\cite{FoKa1972}, random line
representation for Ising~\cite{Ai1982, PfVe1997a,PfVe1999}, or more generally
random walk representation of $N$-vector models~\cite{BrFrSp1982}, etc...
However, it might not always be easy, or even possible, to establish
properties~\eqref{eq:PUB}, 
\eqref{eq:FE}, \eqref{eq:Spl} and~\eqref{eq:ExD} for the
corresponding weights.

Before proceeding  to explain how such an expansion is used in order to prove
Theorem~\ref{theorem}, let us first discuss the similar, but much simpler case of
self-avoiding walks (SAW), which has been treated, though with a
somewhat different approach, in~\cite{ChCh1986,Io1998}.

\vskip 0.1cm
\noindent
{\bf Self-avoiding Walks.}   
Here we are interested in exact asymptotics of the following quantity:
\begin{equation*}
g_\beta^\mathrm{SAW}(x) =
\sum_{\lambda:\, 0\to x} e^{-\beta|\lambda|}\,,
\end{equation*}
where the sum runs 
over all finite-range SAW connecting $0$ and $x$, $|\lambda|$
is the length (i.e. the number of steps) of $\lambda$, and $\beta>0$. It is
 known~\cite{MaSl1993} that $g_\beta^\mathrm{SAW}$ is well-defined and, 
accordingly, that \eqref{eq:PUB} holds  for all
$\beta$ as soon as $\beta>\beta_c^\mathrm{SAW}(d)$.
The remaining three properties trivially follow with $C_2 = \beta$, 
 $C_3 = 1$ and $C_4 = 0$.

We are now going to decompose a path $\lambda:\, 0\to x$ into
elementary irreducible pieces. The notion of irreducibility is adjusted
to the direction of the  target point $x$ through the geometry 
of the inverse correlation length $\xi_\beta^{\mathrm SAW}$: Let 
${\mathbf K}_\beta^{\mathrm SAW}$ be the convex set supported by 
$\xi_\beta^{\mathrm SAW}$. Choose a dual point 
$t\in \partial{\mathbf K}_\beta^{\mathrm SAW}$ such that 
$(t,x)= \xi_\beta^{\mathrm SAW} (x)$. 
We say
that a vertex $i\in\lambda$ is a 
\textit{
regeneration point} of $\lambda$ if the hyperplane
through $i$, orthogonal to $t$, cuts $\lambda$ into two pieces. 
A path is then
said to be \textit{
 irreducible} if it does not contain any
regeneration points. Since SAW-weights possess the factorization 
property $e^{-\beta |\gamma\amalg\lambda|} =
e^{-\beta|\gamma|}\,e^{-\beta |\lambda |}$,
we arrive to the following  Ornstein-Zernike type equation, or in the
probabilistic jargon a \textit{renewal equation}:
\begin{equation}
\label{eq:renewal}
g_\beta^\mathrm{SAW}(x) = \sum_{y\in\mathbb{Z}^d} c_\beta^\mathrm{SAW}(x-y)\,
g_\beta^\mathrm{SAW}(y)\,,
\end{equation}
where the direct two-point function is defined via the summation of
 path weights over {\em irreducible} paths, 
$$
c_\beta^\mathrm{SAW}(y) = \sum_{\substack{
\lambda:\, 0 \to y\\ \mathrm{irreducible}}} e^{-\beta
|\lambda|}\, .
$$ 

The short range nature of $c_\beta^\mathrm{SAW}$ finds then its precise 
mathematical expression in the claim that the  direct correlation  length
is strictly smaller than $1/\xi_\beta^\mathrm{SAW}(x)$ : 
\begin{equation}
\label{eq:mass}
\liminf_{k\to\infty} -(k|x|)^{-1}\, \log c_\beta^\mathrm{SAW}([kx]) >
\xi_\beta^\mathrm{SAW}(x)\,.
\end{equation}

\begin{figure}[t!]
\onefigure{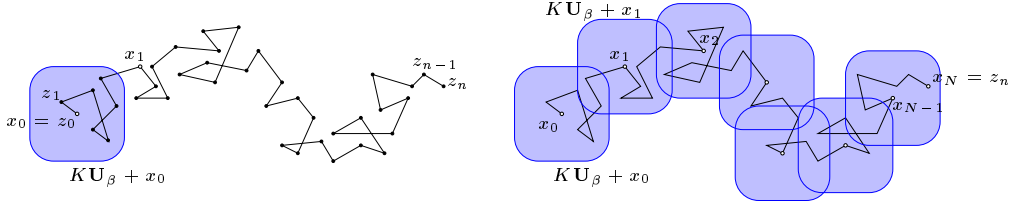}
\caption{The construction of the $K$-skeleton
$\lambda_K=(x_0,\dots,x_N)$ of a path $\lambda=(z_0,\dots,z_n)$. Set $x_0\equiv z_0$; then define iteratively $x_{k+1}$ to be the first point of $\lambda$ outside the set $K\mathbf{U}_\beta + x_k$. When this procedure stops set $x_N\equiv z_n$.}
\label{fig_skeleton}
\end{figure}

It is very easy to establish~\eqref{eq:mass} at large values of $\beta$. 
The main
difficulty is to give a nonperturbative proof. 
As one gets close to $\beta_c^\mathrm{SAW}$, paths
typically have a very complicated messy 
structure at the lattice scale, and one 
expects them to behave properly only at distances large compared to the
correlation length. It is 
thus natural to introduce a coarsegrained description
of the microscopic paths. To this end, we choose some big number $K$, and
construct 
the \textit{$K$-skeleton} $\lambda_K$ of a path $\lambda$ as described
in Fig.~\ref{fig_skeleton}. Notice that $K\mathbf{U}_\beta$ are balls
of radius $K$ in the metric $\xi_\beta^\mathrm{SAW}$, so 
 both the scale and the geometry 
 of the inverse correlation length $\xi_\beta$ literally 
 set up the stage for our  
path coarsegraining procedures.
Given a target point $x$ we use a dual direction
$t\in\partial {\mathbf K}_\beta^{\mathrm SAW} ;\ (t,x) =\xi_\beta (x),$
to measure the 
amount of backtracking by different skeleton steps: For any 
$v\in{\mathbb Z}^d$ define the surcharge cost 
$\mathfrak{s}_t (v) = \xi_\beta^{\mathrm SAW} (v) - (t,v )\geq 0$.
Accordingly, 
the surcharge 
cost of a skeleton $\lambda_K :0\mapsto x$ is defined as sum of all
the surcharge costs of all its steps, 
$\mathfrak{s}_t (\lambda_K) =\sum_k \mathfrak{s}_t (x_k -x_{k-1})$. 

The following energy-entropy argument explains what is gained
by the above change of scale: 
By the splitting property the weight of a skeleton $\lambda_K : 0\mapsto x$
is bounded above as
$$
q_\beta^\mathrm{SAW}(\lambda_K) = 
\sum_{\lambda\sim\lambda_K} e^{-\beta |\lambda|} 
\leq {\mathrm e}^{ -\xi_\beta^{\mathrm SAW} (x) -\mathfrak{s}_t (\lambda_K)} .
$$
On the other hand, the total number of all $K$-skeletons of length
 $N$ 
is of the
order $(K^{d-1})^N = \exp\{c\, N\log K\}$. Therefore, energy dominates
entropy as soon as $K$ is large. In particular,   up to
exponentially  small probabilities, only a negligible fraction of 
$\lambda_K$-steps
 are backtracking (have surcharges of order $K$) with respect to $x$. 
 Because of
 the finite 
energy condition \eqref{eq:FE}
 such predominant forward structure
of paths on the fixed finite scale $K$ implies the massgap property 
\eqref{eq:mass}. 

Repeated iterations of \eqref{eq:renewal} give:
\begin{equation*} \label{eq:RW}
g_\beta^\mathrm{SAW}(x)\, e^{\xi_\beta (x)}= \sum_{n=1}^{\infty}
\sum_{y_1 +y_2 +\dots y_n =x}
\, \prod_{i=1}^n
c_\beta^\mathrm{SAW}(y_i)\, e^{(t,y_i)}\,
\end{equation*}
which, in probabilistic terms, is tantamount to independence
of different irreducible paths.  Since $c_\beta^\mathrm{SAW}(y)\,
e^{(t,y)}$ is a probability measure on $\mathbb{Z}^d$ when
$t\in\partial{\mathbf K}_\beta^{\mathrm SAW}$, the proof of the OZ-formula for
self-avoiding walks reduces to a local limit computation for sums of
independent random variables with exponentially decaying tails. A very similar
construction applies to all subcritical short range Bernoulli bond
percolation models \cite{CaChCh1991,CaIo2002}.  
\vskip 0.2cm

\noindent
{\bf Ising Model.} 
A convenient graphical representation for the pair-correlation functions of
the Ising model is the random-line
representation~\cite{PfVe1997a,PfVe1999,Ai1982},
\begin{equation*}
g_\beta(x) = \langle\sigma_x\,\sigma_y\rangle_\beta = \sum_{\lambda:\,
x\to y} q_\beta(\lambda)\,,
\end{equation*}
where the sum runs over admissible paths connecting $x$ and $y$, making
jumps only between sites $u$ and $v$ such that $J_{u-v}>0$, and $q_\beta$ is
some strictly positive weight.
It can then be proved that properties~\eqref{eq:PUB}, \eqref{eq:FE}, \eqref{eq:Spl}
and~\eqref{eq:ExD} hold for this representation. We would like to proceed
similarly to what was done for SAW. The major new difficulty 
arising now is 
 that
for the Ising
model  the weights do not factorize: $q_\beta(\lambda \amalg\gamma) \neq
q_\beta(\lambda )\, q_\beta(\gamma)$. The closest expression
to~\eqref{eq:renewal} one can write now is
\begin{equation}
\label{eq:pseudorenewal}
g_\beta(x) = \sum_{y\in\mathbb{Z}^d}
\sum_{\substack{\lambda :\, y\to
x\\ \text{irreducible}}} q_\beta(\lambda)
\sum_{\substack{\gamma :\, 0\to y\\ \gamma\amalg\lambda
\text { admissible}}} q_\beta(\gamma \,|\, \lambda )\,.
\end{equation}
The admissibility constraint does not create serious problems, for example it 
 is automatically satisfied across the  regeneration points.  
However, the presence of the conditional 
weight $q_\beta(\gamma \,|\,
\lambda)$ instead of $q_\beta(\gamma)$ destroys the independence and, 
consequently,  dramatically changes  the  probabilistic structure
of the equation.

At this point one is compelled to adjust  the traditional point of view 
on the nature of the OZ-formula and to try to understand it not through the 
spectral theory of the renewal type equations \eqref{eq:OZrenewal} or 
\eqref{eq:renewal}, but rather in a broader context of 
 local limit  properties    
 of one-dimensional systems of objects (irreducible paths) 
under appropriate mixing conditions.  As we shall explain below such an
 approach leads to a reformulation of the problem in terms of the
statistical 
mechanics of Ruelle's operator for shifts over countable 
alphabets. The OZ-formula \eqref{eq:OZising} is recovered in this way from a 
local limit analysis based on the analytic perturbation theory of the 
corresponding Perron-Frobenius eigenvalue. 

For the moment let us postpone the exact definition of path-irreducibility
and focus on the induced irreducible representation
 of $g_\beta (x)$:
\begin{equation}
\label{eq:decomp}
g_\beta(x) = \sum_{n\geq 1} \sum_{\substack{\lambda:\, 0\to x \\
\lambda=\lambda_1\amalg\ldots\amalg\lambda_n\\ \lambda_k \text{ irreducible},\,
k=1,\ldots, n}} q_\beta(\lambda_1\amalg\ldots\amalg\lambda_n)\, .
\end{equation}
Repeated iterations of \eqref{eq:pseudorenewal} suggest to rewrite
the weights as 
\begin{equation*}
 q_\beta(\lambda_1\amalg\ldots\amalg\lambda_n)
 = q_\beta(\lambda_n)\prod_{k=1}^{n-1} q_\beta(\lambda_{k} \,|\,
\lambda_{k+1}\amalg \cdots \amalg\lambda_n)
\end{equation*}
Introducing a dummy empty path $\emptyset$ and using it to extend
finite strings of irreducible paths $(\lambda_1 ,\dots,\lambda_n )$
 to infinite strings $\underline{\lambda} =
(\lambda_1 ,\dots,\lambda_n ,\emptyset , \emptyset ,\dots )$, we rewrite
\eqref{eq:decomp}
in a more  uniform and compact form:
\begin{equation}
\label{eq:ruelle1}
g_\beta(x) = \sum_{n\geq 1} \sum_{\underline{\lambda} :0\mapsto x}
\,\exp \{ \psi (\tau^k\underline{\lambda} )\}\, ,
\end{equation}
where  
$\tau$ is the shift 
$\tau\underline{\lambda} = 
\tau (\lambda_1,\lambda_2,\dots )=(\lambda_2 ,\dots )$
 and 
the potential $\psi$ is defined via 
${\mathrm e}^{\psi (\underline{\lambda} )} =q_\beta (\lambda_1\, |
\lambda_2\amalg\dots )$ under the convention 
$q_\beta (\lambda \, |\, \emptyset\amalg\emptyset\amalg\dots )=
q_\beta (\lambda )$. 

The expansion \eqref{eq:ruelle1} suggests to introduce the normalized 
operator $L$ as follows:  
\begin{equation}
\label{eq:ruelle2}
Lf(\underline{\lambda} 
  ) \, =\,
\sum_{\nu\,\text{irreducible}}{\mathrm e}^{\psi(\nu,\underline{\lambda}
) +
(t,V(\nu )) }
f(\nu ,\underline{\lambda} 
) \, ,
\end{equation}
where $V(\nu )\in {\mathbb Z}^d$ is the displacement along the path $\nu$, 
 and, as before, 
  $t\in \partial {\mathbf K}_\beta$ satisfies $(t,x) = \xi_\beta (x )$.

\vskip 0.1cm
\eqref{eq:ruelle1} falls into the framework of the classical Gaussian local 
limit theory once the potential $\psi$ happens to be H\"{o}lder
continuous: Given two different strings 
$\underline{\lambda} = (\lambda_1 ,\lambda_2,\dots )$ and 
$\underline{\lambda}^\prime = (\lambda_1^\prime ,\lambda_2^\prime,\dots )$ of
 irreducible paths  define ${\mathbf i} ( \underline{\lambda}, 
\underline{\lambda}^\prime) = 
\min\{ k\geq 1 :\, \lambda_k\neq \lambda_k^\prime \} $. We would like to find
$\theta \in (0,1)$ such that 
\begin{equation}
\label{psi_norm}
 |\psi |_\theta \, =\, \inf_{k} \inf_{ {\mathbf i} ( \underline{\lambda},
\underline{\lambda}^\prime)>k} 
\frac{\big| \psi (\underline{\lambda}) - 
\psi (\underline{\lambda^\prime})\big|}{\theta^k} \, <\,\infty .
\end{equation}
Clearly, the exponential mixing property \eqref{eq:ExD} alone is not 
sufficient to ensure \eqref{psi_norm}: Even if $ k\gg 1$ the sum 
$\sum_{u\in \lambda_1}\sum_{v\in \lambda_k} \theta^{|v-u |}$ is beyond
control unless one imposes further restrictions on the geometry 
of irreducible paths. This is precisely the motivation behind the following
 refined definition of irreducible paths: Given $t\in \partial
{\mathbf K}_\beta$ and a renormalization scale $K$, let us say that
$i_l$ is a break point of a path $\lambda = (i_1,\dots ,i_n )$ if
$i_l$ is a regeneration point of $\lambda$ and
$$
 \{i_{l+1},\ldots, i_n\} 
\subset i_l + K\mathrm{U}_\beta + \mathcal{C}_\delta(t)\,,
$$
where $\mathcal{C}_\delta(t) = \{v\in\mathbb{Z}^d \,:\,
( t,v)  > (1-\delta ) \xi_\beta (v)\}$ is a positive cone 
 along the direction of $t$.
A path $\lambda$ is 
\textit{irreducible} if it does not contain break points. 
If $\lambda_1\amalg\cdots \amalg \lambda_n$ is the 
concatenation of such irreducible
paths, then for every $1\leq l < k\leq n$, 
$$
\sum_{u\in \lambda_l}\sum_{v\in \lambda_k} \theta^{|v-u |} \, \leq
\, C_5 (K) \theta^{k-l} , 
$$
which, in view of \eqref{eq:ExD}, already implies that that the 
potential $\psi$ in \eqref{eq:ruelle2} satisfies the H\"{o}lder
 condition
 \eqref{psi_norm}. It remains to check that up to exponentially 
negligible weights typical paths $\gamma :\, 0\mapsto x$ contain 
a density of break points. However, 
 coarsegraining procedures similar to those developed above  for SAW,
 and 
based on the 
properties~\eqref{eq:PUB}, \eqref{eq:FE} and~\eqref{eq:Spl}, show that there exist
$M<\infty$ and $\nu>0$ such that 
\begin{equation*}
\sum_{\substack{\lambda:\, 0\to y\\ \text{irreducible}}} q_\beta(\lambda)
\leq M\, e^{-\nu |x|}\, e^{-\xi_\beta(x)}\,,
 \end{equation*}
This is the
analogue of the mass separation result~\eqref{eq:mass}.

Both the local limit asymptotics of Theorem~\ref{theorem}  
  and the geometry of ${\mathbf K}_\beta$ can be 
now read from the analytic dependence  of the leading eigenvalue 
$\rho (z)$  of $L_z$ on the perturbation $z\in {\mathbb C}^d$ in 
$$
L_z f (\underline{\lambda} )\, = \, L\left( {\mathrm e}^{(z, V(\lambda_1 ))}
f\right)\, .
$$
Indeed, $\log\rho (z)$ is nothing but the limiting log-moment 
generation function of $1/n \sum_{1}^{n} V(\lambda_i )$. 
On the other hand, since ${\mathbf K}_\beta$ is the domain of 
convergence of \eqref{eq:Suscept}, 
the equation of $\partial {\mathbf K}_\beta$ in a neighborhood
of $t$ is given by $\{t +z\in {\mathbb R}^d\, :\, \rho (z) =1\}$. In other 
words, locally the surface $\partial {\mathbf K}_\beta$ is just a level
set of $\rho$ through $t$. Furthermore, choosing perturbations 
$z$ of the form $z= p_1 n_x +i{\mathbf p}$, where ${\mathbf p}$ lies in 
the tangent hyperplane to $\partial {\mathbf K}_\beta$ at $t$, we  
recover the results of Paes-Leme (\cite{Pa1978}, Proposition~4.2) on
the analyticity of the 1-particle mass shells.

\acknowledgments{M.C. and Y.V. gratefully acknowledge the kind
hospitality of Technion where part of this work was done.}

\end{document}